\documentclass[10pt,a4paper]{amsart}
\usepackage{amsfonts,amsmath,amssymb}
\usepackage{hyperref}

\newtheorem*{theorem*}{Theorem}
\newtheorem{lemma}{Lemma}[subsection]
\newtheorem{proposition}[lemma]{Proposition}

\newtheorem{theorem}[lemma]{Theorem}
\newtheorem{definition}[lemma]{Definition}
\newtheorem{notation}[lemma]{Notation}

\newtheorem{corollary}[lemma]{Corollary}

\baselineskip 18pt \textwidth 16cm  \theoremstyle{plain}

\newcommand{\tr}{\operatorname{Tr}}

\newcommand{\cc}{\mathbb{C}}

\newcommand{\eps}{\varepsilon}

\newcommand{\Z}{{\mathbb Z}}

\newcommand{\R}{{\mathbb R}}

\newcommand{\Fre}{{Fr\'{e}chet \,}}

\newcommand{\g}{{\mathfrak{g}}}
\newcommand{\h}{{\mathfrak{h}}}

\newcommand{\gd}{\g^{\sigma}}

\begin{document}

\author{Avraham Aizenbud}
\address{Avraham Aizenbud and Dmitry Gourevitch, Faculty of Mathematics
and Computer Science, The Weizmann Institute of Science POB 26,
Rehovot 76100, ISRAEL.} \email{aizenr@yahoo.com}
\author{Dmitry Gourevitch}
\email{guredim@yahoo.com}
\title[An archimedean analog of Jacquet - Rallis theorem]{An archimedean analog of Jacquet - Rallis theorem}
\begin{abstract} In this paper we prove that the symmetric pair
$(GL_{n+k}(F),GL_{n}(F) \times GL_k(F))$ is a Gelfand pair for any
local field $F$ of characteristic 0. For non-archimedean $F$ it
has been proven in \cite{JR}. We use techniques developed in
\cite{AG2} to generalize their proof to general local fields.
\end{abstract}


%
%
%
%
%
%


\maketitle \setcounter{tocdepth}{1}
 \tableofcontents

\section{Introduction}
Fix a local field $F$ of characteristic zero. Consider the
standard embedding of $GL_{n} \times GL_k$ to $GL_{n+k}$. The goal
of this paper is to prove the following theorem.

\begin{theorem*}[I]
The symmetric pair $(GL_{n+k},GL_{n} \times GL_k)$ is a Gelfand
pair. Namely, for any irreducible admissible representation
$(\pi,E)$ of $\mathrm{GL}_{n+k}(F)$ we have
\begin{equation}\label{dim1}
\dim Hom_{GL_{n}(F) \times GL_k(F)}(E^{\infty},\cc) \leq 1.
\end{equation}
\end{theorem*}

Using the classical arguments of Gelfand and Kazhdan we deduce it
(in section \ref{GelPairs}) from the following theorem.

\begin{theorem*}[II]
Any distribution $\xi$ on $GL_{n+k}(F)$ which is invariant with
respect to two-sided action of $GL_{n}(F) \times GL_k(F)$ is
invariant with respect to the antiinvolution $g \mapsto g^{-1}$.
\end{theorem*}

In the non-archimedean case these theorems had been proven in
\cite{JR}. The techniques developed in \cite{AG2} enabled us to
generalize that proof to work in the archimedean case. Also, those
techniques made the proof shorter by eliminating some of the
technical computations.

In the case $k=1$, already the pair $(GL_n \times \{1\},GL_{n+1})$
is a Gelfand pair (see \cite{AGS}). For non-archimedean $F$ it is
even proven to be a strong Gelfand pair (see \cite{AGRS}).

\subsection*{Structure of the paper}

In subsection 2.1 we give preliminaries on Gelfand pairs and their
connection to invariant distributions. In subsections 2.2-2.4 we
introduce those techniques developed in \cite{AG2} which are
relevant to symmetric pairs. We also formulate theorem II' which
is equivalent to theorem II. In section 3 we prove theorem II'.

\subsection*{Acknowledgements}
We would like to thank our teacher \textbf{Joseph Bernstein} for
our mathematical education. We also thank \textbf{Joseph
Bernstein}, \textbf{Herve Jacquet}, \textbf{Erez Lapid}, and
\textbf{Eitan Sayag} for fruitful discussions.

\section{Preliminaries and notations}
\begin{itemize}
\item All the algebraic varieties and algebraic
groups that we will consider will be defined over $F$.
\item For a group $G$ acting on a set $X$ and an element $x \in X$
we denote by $G_x$ the stabilizer of $x$.
\item By a reductive group we mean an algebraic reductive group.
\end{itemize}

In this paper we will refer to distributions on algebraic
varieties over archimedean and non-archimedean fields. In the
non-archimedean case we mean the notion of distributions on
$l$-spaces from \cite{BZ}, that is linear functionals on the space
of locally constant compactly supported functions. In the
archimedean case one can consider the usual notion of
distributions, that is continuous functionals on the space of
smooth compactly supported functions, or the notion of Schwartz
distributions (see e.g. \cite{AG1}). It does not matter here which
notion to choose since in the cases of consideration of this
paper, if there are no nonzero equivariant Schwartz distributions
then there are no nonzero equivariant distributions at all (see
Theorem 4.0.8 in \cite{AG2}).

\subsection{Gelfand pairs} \label{GelPairs}
$ $\\
In this section we recall a technique due to Gelfand and Kazhdan
(\cite{GK}) which allows to deduce statements in representation
theory from statements on invariant distributions. For more
detailed description see \cite{AGS}, section 2.

\begin{definition}
Let $G$ be a reductive group. By an \textbf{admissible
representation of} $G$ we mean an admissible representation of
$G(F)$ if $F$ is non-archimedean (see \cite{BZ}) and admissible
smooth \Fre representation of $G(F)$ if $F$ is archimedean.
\end{definition}

We now introduce three notions of Gelfand pair.

\begin{definition}\label{GPs}
Let $H \subset G$ be a pair of reductive groups.
\begin{itemize}
\item We say that $(G,H)$ satisfy {\bf GP1} if for any irreducible
admissible representation $(\pi,E)$ of $G$ we have
$$\dim Hom_{H(F)}(E,\cc) \leq 1$$

\item We say that $(G,H)$ satisfy {\bf GP2} if for any irreducible
admissible representation $(\pi,E)$ of $G$ we have
$$\dim Hom_{H(F)}(E,\cc) \cdot \dim Hom_{H}(\widetilde{E},\cc)\leq
1$$

\item We say that $(G,H)$ satisfy {\bf GP3} if for any irreducible
{\bf unitary} representation $(\pi,\mathcal{H})$ of $G(F)$ on a
Hilbert space $\mathcal{H}$ we have
$$\dim Hom_{H(F)}(\mathcal{H}^{\infty},\cc) \leq 1.$$
\end{itemize}

\end{definition}
Property GP1 was established by Gelfand and Kazhdan in certain
$p$-adic cases (see \cite{GK}). Property GP2 was introduced in
\cite{Gross} in the $p$-adic setting. Property GP3 was studied
extensively by various authors under the name {\bf generalized
Gelfand pair} both in the real and $p$-adic settings (see
e.g.\cite{vD,vD-P,Bos-vD}).

We have the following straightforward proposition.

\begin{proposition}
$GP1 \Rightarrow GP2 \Rightarrow GP3.$
\end{proposition}

We will use the following theorem from \cite{AGS} which is a
version of a classical theorem of Gelfand and Kazhdan.

\begin{theorem}\label{DistCrit}
Let $H \subset G$ be reductive groups and let $\tau$ be an
involutive anti-automorphism of $G$ and assume that $\tau(H)=H$.
Suppose $\tau(\xi)=\xi$ for all bi $H(F)$-invariant distributions
$\xi$ on $G(F)$. Then $(G,H)$ satisfies GP2.
\end{theorem}

In our case GP2 is equivalent to GP1 by the following proposition.

\begin{proposition}
Suppose $H \subset \mathrm{GL}_{n}$ is transpose invariant
subgroup. Then $GP1$ is equivalent to $GP2$ for the pair
$(\mathrm{GL}_{n},H)$.
\end{proposition}

For proof see \cite{AGS}, proposition 2.4.1.

\begin{corollary}
Theorem II implies Theorem I.
\end{corollary}

\subsection{Symmetric pairs}
$ $\\
In this subsection we review some tools developed in \cite{AG2}
that enable to prove that a symmetric pair is a Gelfand pair.

\begin{definition}
A \textbf{symmetric pair} is a triple $(G,H,\theta)$ where $H
\subset G$ are reductive groups, and $\theta$ is an involution of
$G$ such that $H = G^{\theta}$. We call a symmetric pair
\textbf{connected} if $G/H$ is connected.

For a symmetric pair $(G,H,\theta)$ we define an antiinvolution
$\sigma :G \to G$ by $\sigma(g):=\theta(g^{-1})$, denote $\g:=Lie
G$, $\h := LieH$, $\gd:=\{a \in \g | \theta(a)=-a\}$. Note that
$H$ acts on $\gd$ by the adjoint action. Denote also
$G^{\sigma}:=\{g \in G| \sigma(g)=g\}$ and define a
\textbf{symmetrization map} $s:G \to G^{\sigma}$ by $s(g):=g
\sigma(g)$.
\end{definition}
\begin{definition}
Let $(G_1,H_1,\theta_1)$ and $(G_2,H_2,\theta_2)$ be symmetric
pairs. We define their \textbf{product} to be the symmetric pair
$(G_1 \times G_2,H_1 \times H_2,\theta_1 \times \theta_2)$.
\end{definition}

\begin{definition}
We call a symmetric pair $(G,H,\theta)$ \textbf{good} if for any
closed $H(F) \times H(F)$ orbit $O \subset G(F)$, we have
$\sigma(O)=O$.
\end{definition}

\begin{proposition} \label{GoodCrit}
Let $(G,H,\theta)$ be a connected symmetric pair. Suppose that for
any $g \in G(F)$ with closed orbit, $H^1(F,(H \times H)_g)$ is
trivial. Then the pair $(G,H,\theta)$ is good.
\end{proposition}
For proof see \cite{AG2}, Corollary 7.1.5.

\begin{definition}
We say that a symmetric pair $(G,H,\theta)$ is a \textbf{GK pair}
if any $H(F) \times H(F)$ - invariant distribution on $G(F)$ is
$\sigma$ - invariant.
\end{definition}
\begin{definition}
We define an involution $\theta_{n,k}:GL_{n+k} \to GL_{n+k}$ by
$\theta_{n,k}(x)=\eps x \eps$ where $\eps = \begin{pmatrix}
  Id_{nn} & 0_{nk} \\
  0_{kn} & -Id_{kk}
\end{pmatrix}$. Note that $(GL_{n+k},GL_{n} \times GL_{k},\theta_{n,k})$ is a symmetric pair.
If there is no ambiguity we will denote $\theta_{n,k}$ just by
$\theta$.
\end{definition}

Theorem II is equivalent to the following theorem.

\begin{theorem*}[II']
The pair $(GL_{n+k},GL_{n}\times GL_{k},\theta_{n,k})$ is a GK
pair.
\end{theorem*}

\subsection{Descendants of symmetric pairs}
\begin{proposition} \label{PropDescend}
Let  $(G,H,\theta)$ be a symmetric pair. Let $g \in G(F)$ such
that $HgH$ is closed. Let $x=s(g)$. Then $x$ is semisimple.
\end{proposition}
For proof see e.g. \cite{AG2}, Proposition 7.2.1.
\begin{definition}
In the notations of the previous proposition we will say that the
pair $(G_x,H_x,\theta|_{G_x})$ is a descendant of $(G,H,\theta)$.
\end{definition}

\subsection{Regular symmetric pairs}

\begin{notation}
Let $V$ be an algebraic finite dimensional representation over $F$
of a reductive group $G$. Denote $Q(V):=V/V^G$. Since $G$ is
reductive, there is a canonical embedding $Q(V) \hookrightarrow
V$.
\end{notation}

\begin{notation}
Let $(G,H,\theta)$ be a symmetric pair. We denote by
$\mathcal{N}_{G,H}$ the subset of all the nilpotent elements in
$Q(\gd)$. Denote $R_{G,H}:=Q(\gd) - \mathcal{N}_{G,H}$.
\end{notation}
Our notion of $R_{G,H}$ coincides with the notion $R(\gd)$ used in
\cite{AG2}, Notation 2.1.10. This follows from Lemma 7.1.11 in
\cite{AG2}.
\begin{definition}
Let $\pi$ be an action of a reductive group $G$ on a smooth affine
variety $X$.
We say that an algebraic automorphism $\tau$ of $X$ is \textbf{$G$-admissible} if \\
(i) $\pi(G(F))$ is of index at most 2 in the group of
automorphisms of $X$
generated by $\pi(G(F))$ and $\tau$.\\
(ii) For any closed $G(F)$ orbit $O \subset X(F)$, we have
$\tau(O)=O$.
\end{definition}

\begin{definition}
Let $(G,H,\theta)$ be a symmetric pair. We call an element $g \in
G(F)$ \textbf{admissible} if\\
(i) $Ad(g)$ commutes with $\theta$ (or, equivalently, $s(g)\in Z(G)$) and \\
(ii) $Ad(g)|_{\g^{\sigma}}$ is $H$-admissible.
\end{definition}

\begin{definition}
We call a symmetric pair $(G,H,\theta)$ \textbf{regular} if for
any admissible $g \in G(F)$ such that every $H(F)$-invariant
distribution on $R_{G,H}$ is also $Ad(g)$-invariant,
we have\\
(*) every $H(F)$-invariant distribution on $Q(\gd)$ is also
$Ad(g)$-invariant.
\end{definition}

Clearly product of regular pairs is regular (see \cite{AG2},
Proposition 7.4.4).

\begin{theorem} \label{GoodHerRegGK}
Let $(G,H,\theta)$ be a good symmetric pair such that all its
descendants (including itself) are regular. Then it is a GK pair.
\end{theorem}
For proof see \cite{AG2}, Theorem 7.4.5.

Now we would like to formulate a regularity criterion. For it we
will need the following lemma and notation.

\begin{lemma}
Let $(G,H,\theta)$ be a symmetric pair. Then any nilpotent element
$x \in \gd$ can be extended to an $sl_2$ triple $(x,d(x),x_-)$
such that $d(x) \in \h$ and $x_- \in \gd$.
\end{lemma}
For proof see e.g. \cite{AG2}, Lemma 7.1.11.

\begin{notation}
We will use the notation $d(x)$ from the last lemma in the future.
It is not uniquely defined but whenever we will use this notation
nothing will depend on its choice.
\end{notation}

\begin{proposition} \label{SpecCrit}
Let $(G,H,\theta)$ be a symmetric pair. Suppose that for any
nilpotent $x \in \gd$ we have
$$\tr(ad(d(x))|_{\h_x}) < \dim \gd .$$
Then the pair $(G,H,\theta)$ is regular.
\end{proposition}
This proposition follows from \cite{AG2} (Propositions 7.3.7 and
7.3.5 and Remark 7.4.3).

\begin{notation}
Let $E$ be a quadratic extension of $F$. Let $G$ be an algebraic
group defined over $F$. We denote by $G_{E/F}$ the canonical
algebraic group defined over $F$ such that $G_{E/F}(F)=G(E)$.
\end{notation}

\begin{theorem}\label{2RegPairs}
Let $G$ be a reductive group.\\
(i)Consider the involution $\theta$ of $G\times G$ given by
$\theta((g,h)):= (h,g)$. Its fixed points form the diagonal
subgroup $\Delta G$. Then the symmetric pair $(G \times G, \Delta
G, \theta)$ is
regular.\\
(ii) Let $E$ be a quadratic extension of $F$. Consider the
involution $\gamma$ of $G_{E/F}$ given by the nontrivial element
of $Gal(E/F)$. Its fixed points form $G$. Then the symmetric pair
$(G_{E/F}, G, \gamma)$ is regular.
\end{theorem}

This theorem follows from \cite{AG2} (Theorem 7.6.5, Proposition
7.3.5 and Remark 7.4.3).

\section{Proof of theorem II'}
By theorem \ref{GoodHerRegGK} it is enough to prove that our pair
is good and all its descendants are regular.

In the first subsection we will compute the descendants of our
pair. Proposition \ref{GoodCrit} will imply that our pair is good.

In the second subsection we will prove that all the descendants
are regular using theorem \ref{2RegPairs} and proposition
\ref{SpecCrit}.
\subsection{Computation of descendants of the pair $(GL_{n+k},GL_{n} \times GL_k)$}

\begin{theorem} \label{ComputeDescent}
All the descendants of the pair $(GL_{n+k},GL_{n}\times
GL_{k},\theta_{n,k})$ are products of pairs of the types

(i) $((GL_m)_{E/F} \times (GL_m)_{E/F}, \Delta (GL_m)_{E/F},
\theta)$ for some field extension $E/F$

(ii) $((GL_m)_{E/F}, (GL_m)_{L/F}, \gamma)$ for some field
extension $L/F$ and its quadratic extension $E/L$

(iii) $(GL_{m+l},GL_{m}\times GL_{l},\theta_{m,l}).$
\end{theorem}
\begin{proof}
Let $x \in GL_{n+k}^{\sigma}(F)$ be a semisimple element. We have
to compute $G_x$ and $H_x$. Since $x \in G^{\sigma}$, we have
$\eps x \eps = x^{-1}$. Let $V = F^{n+k}$. Decompose $V :=
\bigoplus _{i=1}^s V_i$ such that the minimal polynomial of
$x|_{V_i}$ is irreducible. Now $G_x(F)$ decomposes to a product of
$GL_{E_i}(V_i)$, where $E_i$ is the extension of $F$ defined by
the minimal polynomial of $x|_{V_i}$ and the $E_i$-vector space
structure on $V_i$ is given by $x$.

Clearly for any $i$, $\eps(V_i) = V_j$ for some $j$.
Now we see that $V$ is a direct sum of spaces of the following two types\\
A. $W_1 \oplus W_2$ such that the minimal polynomials of
$x|_{W_i}$ are irreducible and $\eps(W_1)=W_2$.\\
B. $W$ such that the minimal polynomial of $x|_{W}$ is irreducible
and $\eps(W)=W$.

It is easy to see that in case A we get the symmetric pair (i).

In case B there are two possibilities: (1) $x = x^{-1}$ and (2) $x
\neq x^{-1}$. It is easy to see that in case (1) we get the
symmetric pair (iii) and in case (2) we get the symmetric pair
(ii).
\end{proof}

\begin{corollary}
The pair ($GL_{n+k},GL_{n} \times GL_k$) is good.
\end{corollary}

\subsection{All the descendants of the pair $(GL_{n+k},GL_{n} \times GL_k)$ are regular}
$ $\\
Clearly if a pair $(G,H,\theta)$ is regular then the pair
$(G_{E/F},H_{E/F},\theta)$ is also regular for any field extension
$E/F$. Therefore by Theorem \ref{ComputeDescent} and Theorem
\ref{2RegPairs} it is enough to prove that the pair
$(GL_{n+k},GL_{n}\times GL_{k},\theta_{n,k})$ is regular.

In case $n \neq k $ it follows from the definition since any
admissible $g\in GL_{n+k}$ lies in $GL_n \times GL_k$.

So we can assume $n=k>0$. Hence by proposition \ref{SpecCrit} it
is enough to prove the following key lemma.

\begin{lemma}[Key lemma]\footnote{This lemma is similar to Lemma 3.1 in \cite{JR}, section 3.2. The
proofs are also similar.} \label{Key} Let $x \in
gl_{2n}^{\sigma}(F)$ be a nilpotent element and $d:=d(x)$. Then
$$\tr(ad(d)|_{(gl_{n}(F) \times gl_n(F))_x}) < 2n^2.$$
\end{lemma}

We will need the following definition and lemmas.

\begin{definition}
We fix a grading on $sl_2(F)$ given by $h \in sl_2(F)_0$ and $e,f
\in sl_2(F)_1$ where $(e,h,f)$ is the standard $sl_2$-triple. A
\textbf{graded representation of $sl_2$} is a representation of
$sl_2$ on a graded vector space $V=V_0 \oplus V_1$ such that
$sl_2(F)_i(V_j) \subset V_{i+j}$ where $i,j \in \Z/2\Z$.
\end{definition}

The following lemma is standard.
\begin{lemma}$ $\\
(i) Every irreducible graded representation of $sl_2$ is
irreducible (-as a usual representation of $sl_2$).\\
(ii) Every irreducible representation $V$ of $sl_2$ admits exactly
two gradings. In one grading the highest weight vector lies in
$V_0$ and in the other grading it lies in $V_1$.
\end{lemma}
\begin{lemma} \footnote{This lemma is similar to Lemma 3.2 in \cite{JR} but computes a
different quantity.} Let $V_1$ and $V_2$ be two irreducible graded
representations of an $sl_2$-triple $(e,h,f)$. Let $r_i:=dimV_i$
and let $w_i \in \{-1,+1\}$ be the parity of the highest weight
vectors of $V_i$. Consider $Hom((V_1,V_2)^e)_0$ - the even part of
the space of $e$-equivariant linear maps $V_1 \to V_2$. Let  $$m:=
\tr(h|_{(Hom(V_1,V_2)^e)_0}) + \tr(h|_{(Hom(V_2,V_1)^e)_0}) -
r_1r_2.$$ Then
$$m=\left\{%
\begin{array}{lll}
    -\min(r_1,r_2), & r_1 \neq r_2 & \mod 2; \\
    -\min(r_1,r_2)(1+w_1w_2), & r_1 \equiv r_2 \equiv 0 & \mod 2; \\
    -\min(r_1,r_2) +w_1w_2(\max(r_1,r_2)-1), & r_1 \equiv r_2 \equiv 1 & \mod 2;
\end{array}%
\right.$$
\end{lemma}
This lemma follows by a direct computation from the following
straightforward lemma.
\begin{lemma}
Let $V_{\lambda}^w$ be the irreducible graded representation of
$sl_2$ with highest weight $\lambda$ and highest weight vector of
parity $w\in \{- 1,+1\}$. Then \setcounter{equation}{0}
\begin{align}
& \tr(h|_{(V_{\lambda}^w)^e)_0}) = \lambda\frac{1+w}{2}\\
& (V_{\lambda}^w)^* = V_{\lambda}^{w(-1)^{\lambda}}\\
& V_{\lambda_1}^{w_1} \otimes V_{\lambda_2}^{w_2} =
\bigoplus_{i=0}^{\min(\lambda_1, \lambda_2)}
V_{\lambda_1+\lambda_2 - 2i}^{w_1w_2(-1)^i}.
\end{align}
\end{lemma}

\begin{proof}[Proof of the key lemma]
Let $V_0:=V_1:=F^n$. Let $V:=V_0 \oplus V_1$ be a $\Z/{2\Z}$
graded vector space. We consider $gl_{2n}(F)$ as the $\Z/{2\Z}$
graded Lie algebra $End(V)$. Note that $gl_n(F) \times gl_n(F)$ is
the even part of $End(V)$ with respect to this grading. Consider
$V$ as a graded representation of the $sl_2$ triple $(x,
d,x_{-})$. Decompose $V$ to graded irreducible representations
$W_i$. Let $r_i:=dimW_i$ and $w_i = \pm 1$ be the parity of the
highest weight vector of $W_i$. Note that if $r_i$ is even then
$\dim (W_i \cap V_0)=\dim (W_i \cap V_1)$. If $r_i$ is odd then
$\dim (W_i \cap V_0)=\dim (W_i \cap V_1) +w_i$. Since $\dim V_0 =
\dim V_1$, we get that the number of indexes $i$ such that $r_i$
is odd and $w_i=1$ is equal to the number of indexes $i$ such that
$r_i$ is odd and $w_i=-1$. We denote this number by $l$. Now
$$\tr(ad(d)|_{(gl_{n}(F) \times gl_n(F))_x}) - 2n^2 =
\tr(d|_{(Hom(V,V)^x)_0})-2n^2= 1/2\sum_{i,j} m_{ij},$$
where
$$m_{ij}:=\tr(d|_{(Hom(W_i,W_j)^x)_0}) + \tr(d|_{(Hom(W_j,W_i)^x)_0}) -
r_ir_j.$$ The $m_{ij}$ can be computed using the previous lemma.

As we see from the lemma, if $r_i$ or $r_j$ are even then $m_{ij}$
is non-positive and $m_{ii}$ is negative. Therefore, if all $r_i$
are even then we are done. Otherwise $l>0$ and we can assume that
all $r_i$ are odd. Reorder the spaces $W_i$ such that $w_i = 1$
for $i \leq l$ and $w_i=-1$ for $i>l$. Now

\begin{multline*}
\sum_{i,j}m_{ij} = \sum_{i \leq l, j \leq l} (|r_i-r_j|-1) +
\sum_{i
> l, j > l} (|r_i-r_j|-1) - \sum_{i \leq l, j > l} (r_i+r_j-1) - \sum_{i > l, j \leq l}
(r_i+r_j-1)=\\
=\sum_{i \leq l, j \leq l} |r_i-r_j| + \sum_{i
> l, j > l} |r_i-r_j| - \sum_{i \leq l, j > l} (r_i+r_j) - \sum_{i > l, j \leq l}
(r_i+r_j) <\\
<\sum_{i \leq l, j \leq l} (r_i+r_j) + \sum_{i
> l, j > l} (r_i+r_j) - \sum_{i \leq l, j > l} (r_i+r_j) - \sum_{i > l, j \leq l}
(r_i+r_j)=0.
\end{multline*}

\end{proof}

\end{document}